\documentclass[11pt]{amsart}
\usepackage{eucal}

\theoremstyle{plain}
\newtheorem*{theorem}{Theorem}
\newtheorem*{lemma}{Lemma}
\newtheorem{step}{Step}

\theoremstyle{remark}
\newtheorem*{remark}{Remark}

\newcommand{\C}{\mathbb{C}}

\newcommand{\cS}{\mathcal{S}}
\newcommand{\cT}{\mathcal{T}}
\newcommand{\X}{\tilde{X}}
\DeclareMathOperator{\diam}{diam}

\begin{document}
\title{Diameter preserving linear bijections of $C(X)$}
\author{M\'at\'e Gy\H{o}ry}
\address{Institute of Mathematics\\
         Lajos Kossuth University\\
         4010 Debrecen, P.O. Box 12, Hungary}
\email{gyorym@math.klte.hu}
\author{Lajos Moln\'ar}
\address{Institute of Mathematics\\
         Lajos Kossuth University\\
         4010 Debrecen, P.O. Box 12, Hungary}
\email{molnarl@math.klte.hu}
\thanks{This paper was written when the second author,
        holding a scholarship of the Volkswagen-Stiftung
        of the Konferenz der Deutschen Akademien der Wissenschaften,
        was a visitor at the University of Paderborn, Germany.
        He is very grateful to
        Prof. K.-H.~Indlekofer for his kind hospitality.
        This research was partially supported also by a grant
        from the Ministry of Education, Hungary, Reg. No. FKFP 0304/1997}
\date{July 3, 1997}

\begin{abstract}
The aim of this paper is to solve a linear preserver problem on
the function algebra $C(X)$. We show that in case $X$ is a first
countable compact Hausdorff space, every linear bijection $\phi :C(X)\to
C(X)$ having the property that $\diam(\phi(f)(X))=\diam(f(X))$
$(f\in C(X))$ is of the form

\centerline{$\phi(f)=\tau \cdot f\circ \varphi +t(f)1 \qquad (f\in
C(X))$}

\noindent
where $\tau \in \C$, $|\tau|=1$, $\varphi:X\to X$ is a
homeomorphism and $t:C(X)\to \C$ is a linear functional.
\end{abstract}
\maketitle

Linear preserver problems concern the question of determing all
linear maps on an algebra which leave a given set, function or
relation defined on the underlying algebra invariant. The study of
linear preserver problems on matrix algebras represents one of the most
active research areas in matrix theory (see the survey
paper \cite{LiTs}). In
the last decade a considerable attention has been paid to the infinite
dimensional case as well and the investigations
have resulted in several important results (e.g. \cite{BrSe}).

Concerning function algebras, the main linear preserver problems
studied so far are the characterizations of linear bijections
preserving some given norm, respectively disjointness of the
support (these latter maps are also called separating).
The linear bijections of $C(X)$ preserving the sup-norm are
determined in the famous Banach-Stone theorem.
To mention some recent papers, we
refer to \cite{Wan}, \cite{Wea}, \cite{FoHe} and \cite{HBN}.

One way of measuring a function $f\in C(X)$ is to consider its sup-norm.
One of the other several possibilities that have sense is to measure
the function in question by some data which reflects how large its
range is. For example, this can be done by considering the diameter of
the range. The aim of this paper is to determine all linear
bijections of $C(X)$ which preserve the seminorm $f\mapsto \diam(f(X))$.
For short, we call these maps diameter preserving.
Obviously, every automorphism $f\mapsto f\circ\varphi$ ($\varphi$ is a
homeomorphism) of $C(X)$ is diameter preserving. The diameter of sets in
$\C$ is clearly invariant under rotation and translation. This gives us
that every linear map of the form
\[
\phi (f)=\tau\cdot f\circ \varphi+t(f)1 \qquad (f\in C(X)),
\]
where $\tau \in \C$ is of modulus 1 and $t:C(X)\to \C$ is a
linear
functional, preserves the diameter of the ranges of functions in
$C(X)$.
Our result which follows shows that under a mild condition on the
topological space $X$, every bijective diameter preserving linear map on
$C(X)$ is of the above form.

\begin{theorem}
Let $X$ be a first countable compact Hausdorff space.
A bijective linear map $\phi: C(X)\to C(X)$
is diameter preserving if and only if there exist
a complex number $\tau$ of modulus 1,
a homeomorphism $\varphi :X\to X$ and a linear functional
$t:C(X)\to \C$ with $t(1)\neq -\tau$ so that $\phi$ is of the
form
\begin{equation}\label{E:form}
\phi (f)=\tau \cdot f\circ \varphi+t(f)1 \qquad (f\in C(X)).
\end{equation}
\end{theorem}

\begin{remark}
The statement of our theorem holds true also for the algebra of all
continuous real valued functions on $X$. Then we, of course, have $\tau=
\pm 1$ and $t:C(X)\to \mathbb{R}$. In fact,
in this case the proof can be carried out in a more simple fashion.

We note that, as a particular case, our result gives the form of all
nonsingular linear transformations on $\mathbb{K}^n$ ($\mathbb{K}$
stands for $\mathbb{C}$ or $\mathbb{R}$) which preserve the maximal
distance between the coordinates.
\end{remark}

\begin{proof}[Proof of Theorem]
It is straightforward to check that every linear map $\phi$ of the
form \eqref{E:form} with $\tau, \varphi, t$ being as
in the statement above, is a diameter preserving linear bijection of
$C(X)$.

Now, suppose that $\varphi :C(X) \to C(X)$ is a linear bijection
which preserves the diameter of the ranges of functions in $C(X)$.
Since the theorem is quite easy to verify in the cases when
$X$ has one, respectively two points, in the proof of the necessity
part we suppose that $X$ has at least three points.

Introduce the following notation.
Let $\X$ stand for the collection of all subsets of $X$ having
exactly two elements. For any $f\in C(X)$ let
\begin{gather*}
S(f) =\{ \{x,y\}\in \X\, :\, \vert f(x)-f(y)\vert=
       \diam(f(X)) \}\\
P(f) =\{ (x,y)\in X\times X\, :\,
       \vert f(x)-f(y)\vert= \diam(f(X))\}\\
T(f) =\{ (x,y,u)\in X\times X\times \C \, :\,
       \vert f(x)-f(y)\vert= \diam(f(X)),\, \\
       \hskip 220pt u=f(x)-f(y)\}, \\
\intertext{and for every $\{x,y\}\in \X$ and $u\in \C$ let}
\cS(\{x,y\})    =\{ f\in C(X)\, :\, \{x,y\}\in S(f)\}\\
\cS_s(\{x,y\})  =\{ f\in C(X)\, :\, \{\{x,y\}\}=S(f)\}\\
\cT(x,y,u)      =\{ f\in C(X)\, :\, (x,y,u)\in T(f)\}\\
\cT_s(x,y,u)    =\{ f\in C(X)\, :\, \{(x,y,u), (y,x,-u)\} =T(f)\}. \\
\intertext{Finally, we define}
G(\{x,y\})      =\bigcap \{ S(\phi(f))\, : f\in C(X),\, \{x,y\} \in
                 S(f)\}\\
H(x,y,u)        =\bigcap \{ T(\phi(f)) \, : f\in C(X),\, (x,y,u)\in
                 T(f)\}.
\end{gather*}

Obviously, for every nonconstant function $f\in C(X)$, the sets $S(f)$,
$P(f)$ and
$T(f)$ are nonempty. Since $X$ is first countable, using Uryson's lemma
it is easy to see that to
every pair $x,y\in X, x\neq y$ of points there exists a continuous real
valued function
$f$ from $X$ into $[0,1]$ such that $f(x)=1, f(y)=0$ and $0<f(z)<1$ $(
z\in X, z\neq x, z\neq y)$. This shows that for any element $\{ x, y\}
\in \X$ and $0\neq u\in \C$, the sets $\cS_s(\{ x,y\})$,
$\cT_s(x,y,u)$ are also nonempty.

After this preparation we can begin the proof of the necessity part of
our statement which will be carried out through a series of steps.
The following observation will be used repeatedly in what
follows.

\begin{lemma}
For an arbitrary integer $n\in{\mathbb N}$ and functions
$f_1,\dots ,f_n\in C(X)$ we have
\[
\diam((f_1+\ldots +f_n)(X))=\diam(f_1(X))+\ldots +\diam(f_n(X))
\]
if and only if there exist an
$\{x,y\}\in \X$ and a complex number $v$ of modulus 1 so that
$f_i\in \cT(x,y,\lambda_i v)$ holds true for every $i=1,\dots ,n$,
where $\lambda_i =\diam(f_i(X))$ $(i=1, \ldots, n)$.
\end{lemma}

Let us assume that
\[
\diam((f_1+\ldots +f_n)(X))=\diam(f_1(X))+\ldots+\diam(f_n(X)).
\]
Then there exist $\{x,y\}\in \X$ and a complex number $v$
of modulus $1$ such that
\[
f_1+\ldots +f_n\in \cT(x,y,(\lambda_1+\ldots +\lambda_n) v).
\]
Now, we compute
\begin{gather*}
\lambda_1+\ldots +\lambda_n=
\vert (f_1+\ldots +f_n)(x)-(f_1+\ldots +f_n)(y)\vert \leq\\
\leq \vert f_1(x)-f_1(y)\vert+\ldots +\vert f_n(x)-f_n(y)\vert\leq
\lambda_1 +\ldots +\lambda_n.
\end{gather*}
It readily follows that
$f_i\in \cT(x,y,\lambda_i v)$ $(i=1,\dots ,n)$.
The converse statement of the lemma is trivial.

\begin{step}\label{Sa}
For arbitrary $\{x,y\}\in \X$ and $0\neq u\in \C$ we have
$G(\{x,y\})\neq \emptyset$ and
$H(x,y,u)\neq \emptyset$.
\end{step}

Let $\{x,y\}\in \X$. We first show that
\begin{equation}\label{E:Sa1}
\bigcap \{ P(\phi(f))\, :\, f\in C(X), \{x,y\}\in S(f)\}
\neq\emptyset.
\end{equation}
Since the collection of sets $P(\phi(f))$ in \eqref{E:Sa1}
consists of closed subsets of the compact Hausdorff space $X\times X$,
in order to verify \eqref{E:Sa1}, it is sufficient to show that this
collection has the finite intersection property.
So, let $f_1,\dots ,f_n\in C(X)$ for which $\{x,y\}\in S(f_1),\ldots,
S(f_n)$. Define $u_i=f_i(x)-f_i(y)$ $(i=1,\ldots, n)$.
Then there exist complex numbers $\mu_i$ with $|\mu_i|=1$ for which
$\mu_i u_i\geq 0$.
Since the diameter of the range of $\mu_i f_i$ is $\mu_i
u_i$ and
\begin{gather*}
\vert (\mu_1 f_1+\ldots +\mu_n f_n)(x)-
(\mu_1 f_1+\ldots +\mu_n f_n)(y)\vert=\\
\vert \mu_1(f_1(x)-f_1(y))+\ldots +
\mu_n(f_n(x)-f_n(y))\vert=\\
\vert \mu_1 u_1+\ldots +\mu_n u_n\vert=
\mu_1 u_1+\ldots +\mu_n u_n=
\vert u_1\vert +\ldots +\vert u_n\vert,
\end{gather*}
thus the diameter of the range of $\mu_1 f_1+\ldots
+\mu_n f_n$ is $\vert u_1\vert +\ldots +\vert u_n\vert$.
By the diameter preserving property of $\phi$ it follows that
the diameter of the range of $\phi(\mu_1 f_1+\ldots +\mu_n f_n)$ is
$\vert u_1\vert +\ldots +\vert u_n\vert$ which equals the sum of the
diameters of the ranges of $\phi(\mu_1 f_1),\ldots, \phi(\mu_n f_n)$.
By our lemma we obtain that there exist $\{x_0,y_0\}\in \X$ and a
complex number $v$ of modulus $1$ for which
$\phi (\mu_i f_i)\in \cT(x_0,y_0,\vert u_i\vert v)$
$(i=1,\dots ,n)$. Obviously, we have
$(x_0,y_0)\in P(\phi(\mu_if_i))=P(\phi(f_i))$.
This shows the desired finite intersection property, hence we have
\eqref{E:Sa1}.
Since there is a nonconstant function in $\cS(\{ x,y\})$, every element
of the intersection appearing in \eqref{E:Sa1} has different
coordinates. This gives us that $G(\{ x,y\})\neq \emptyset$.

Let us now prove the remaining assertion \begin{equation}\label{E:Sa2}
H(x,y,u)=\bigcap \{ T(\phi(f))\, :\, f\in C(X),\, (x,y,u)\in T(f)\}
\neq\emptyset.
\end{equation}
It is easy to see that the sets $T(\phi(f))$ are compact subsets of
$X\times X\times \C$. Therefore, just as above, in order to verify
\eqref{E:Sa2}, it is sufficient to check that the system
$\{T(\phi(f)) \, :\, f\in C(X), \, (x,y,u) \in T(f))\}$ has the finite
intersection property.
Let $f_1,\ldots ,f_n\in C(X)$ be such that $(x,y,u)\in
T(f_1),\ldots,
T(f_n)$. We apparently have $ f_1+\ldots +f_n\in \cT(x,y,nu)$ and
hence it follows that $\diam((f_1+\ldots +f_n)(X))=n\vert u\vert$.
From the diameter preserving property of $\phi$ we deduce
\[
\diam(\phi(f_1+\ldots +f_n)(X))=
\diam(\phi(f_1)(X))+\ldots +\diam (\phi(f_n)(X)).
\]
By Lemma, there exist
$\{x_0,y_0\}\in \X$ and a complex number $v$ of modulus 1 such that
$\phi(f_i)\in \cT(x_0,y_0,\vert u\vert v)$ $(i=1,\dots ,n)$.
Plainly, this can be reformulated as $(x_0,y_0,\vert u\vert v)\in
\cap_{i=1}^n T(\phi(f_i))$ verifying the claimed finite
intersection property. Thus, we obtain \eqref{E:Sa2}.

\begin{step}\label{Sb}
If $\{ x_1,y_1\},\{ x_2,y_2\}\in \X$ and
$\{ x_1,y_1\}\neq \{ x_2,y_2\}$, then we have
$G(\{x_1,y_1\})\cap G(\{x_2,y_2\})=\emptyset$.
\end{step}

Let us suppose on the contrary that
$G(\{x_1,y_1\})\cap G(\{x_2,y_2\})\neq \emptyset$.
Clearly, we may assume that $x_1\neq x_2$.
Let $\{x,y\}\in G(\{x_1,y_1\})\cap G(\{x_2,y_2\})$.
By Uryson's lemma there are functions
$f_1\in \cT(x_1,y_1,1)$ and $f_2\in \cT(x_2,y_2,1)$
with disjoint supports and with ranges in $[0,1]$. Let
\[
u_1=\phi(f_1)(x)-\phi(f_1)(y)
\quad \text{and}  \quad
u_2=\phi(f_2)(x)-\phi(f_2)(y).
\]
Then, by the definition of $G$, we have
$\vert u_1\vert=\vert u_2\vert =1$ and
$(x,y,u_1)\in T(\phi(f_1))$, $(x,y,u_2)\in T(\phi(f_2))$.
Let $t\in [-\pi/3,\pi/3]$ be arbitrary and set
$\mu_t=e^{it}$.
We obtain that $f_1+\mu_t f_2\in \cT(x_1,y_1,1)$
and hence $\{x,y\}\in
S(\phi (f_1+\mu_t f_2))$. We can compute
\begin{equation}\label{E:Sb1}
\begin{split}
\vert u_1+\mu_t u_2\vert = &
\vert (\phi(f_1)(x)-\phi(f_1)(y))+ (\phi (\mu_t f_2)(x)-
\phi (\mu_t f_2)(y))\vert=\\
& \vert \phi (f_1+\mu_t f_2)(x)-\phi (f_1+\mu_t f_2)(y)
\vert=1.
\end{split}
\end{equation}
Since $|u_1|=|u_2|=1$ and \eqref{E:Sb1} holds true for every $t\in
[-\pi/3,\pi/3]$, we easily arrive at a contradiction.

\begin{step}\label{Sc}
For every $\{x,y\}\in \X$, the set $G(\{x,y\})$ has
exactly one element in $\X$. The function $G':\X \to \X$
defined by $\{ G'(\{ x,y\})\}=G(\{ x,y\})$ is a bijection.
\end{step}

Let $\{x,y\}\in \X$. Since $G(\{x,y\})$ is nonempty, we can choose
$\{x_0,y_0\}\in \X$ such that \begin{equation}\label{E:Sc1}
\{x_0,y_0\}\in G(\{x,y\}).
\end{equation}
Using the surjectivity of $\phi$, there exists
a function $f\in C(X)$ for which $\phi(f)\in \cT_s(x_0,y_0,1)$.
Let $ \{x_1,y_1\}\in S(f)$. We have
$G(\{x_1,y_1\})\subseteq S(\phi(f))=\{\{x_0,y_0\}\}$.
Therefore, by Step~\ref{Sa} we obtain
$G(\{x_1,y_1\})=\{\{x_0,y_0\}\}$. Using \eqref{E:Sc1} we deduce
$G(\{x_1,y_1\})\subseteq G(\{x,y\})$. By Step~\ref{Sb} this implies
$\{x,y\}=\{x_1,y_1\}$ and we conclude
$G(\{x,y\})=G(\{x_1,y_1\})=\{\{x_0,y_0\}\}$.

Let us prove now that the function $G'$ is bijective.
In view of Step~\ref{Sb} the injectivity is obvious. To see the
surjectivity, let $\{x,y\}\in \X$ and pick an $f\in C(X)$ for which
$\phi(f)\in \cT_s(x,y,1)$. If $\{x_0,y_0\}\in S(f)$, then we have
$G'(\{ x_0,y_0\})\in S(\phi(f))$. It follows that
$G'(\{ x_0,y_0\})=\{ x,y\}$ verifying our claim.

\begin{step}\label{Sd}
For arbitrary $\{x,y\}\in \X$ and $f\in C(X)$, if
$\phi(f)\in \cS_s(G'(\{x,y\}))$, then we have
$f\in \cS_s(\{x,y\})$.
\end{step}

If $\{x_0,y_0\}\in S(f)$ is arbitrary, then
$G(\{x_0,y_0\}) \subseteq S(\phi(f))=G(\{x,y\})$.
By Step~\ref{Sb} it follows that
$\{ x_0,y_0\}=\{ x,y\}$. Consequently, $S(f)=\{ \{x,y\}\}$.

\begin{step}\label{Se}
If $\{x_1,y_1\},\{x_2,y_2\}\in \X$ and
$\{x_1,y_1\}\cap \{x_2,y_2\}\neq\emptyset$, then we have
\[
G'(\{x_1,y_1\})\cap G'(\{x_2,y_2\})\neq\emptyset.
\]
By Step~\ref{Sc} it follows that if $\{x_1,y_1\},\{x_2,y_2\}\in \X$
have exactly one element in common, then the same holds true
for $G'(\{x_1,y_1\})$ and $G'(\{x_2,y_2\})$.
\end{step}

Let $\{x,y_1\},\{x,y_2\}\in \X$ with $y_1\not=y_2$.
Assume that
\[
G'(\{x,y_1\})\cap G'(\{x,y_2\})=\emptyset.
\]
The surjectivity of $\phi$ implies that there exist
functions $\phi(f_1),\phi(f_2)\in C(X)$ with the following
properties. The supports of $\phi(f_1)$ and $\phi(f_2)$ are
disjoint,
\[
\phi(f_1)\in \cS_s(G'(\{x,y_1\})) \quad \text{and} \quad
\phi(f_2)\in \cS_s(G'(\{x,y_2\})),
\]
the ranges of $\phi(f_1), \phi(f_2)$ are included in $[-1/2,1/2]$ and,
finally,
\[
\diam(\phi(f_1)(X))=\diam(\phi(f_2)(X))=1.
\]
By Step~\ref{Sd}, we infer
$f_1\in \cS_s(\{x,y_1\}), f_2\in \cS_s(\{x,y_2\})$.
For every complex number $\mu$ with
$\vert\mu\vert=1$, the diameter of the range of
$\phi(f_1)+\mu\phi(f_2)$ is $1$ and hence the same must hold
true for $f_1+\mu f_2$.

Since $f_2\in \cS_s(\{x,y_2\})$, we have $f_2(y_1)\neq f_2(x)$.
Define
\[
\mu=\frac{f_1(x)-f_1(y_1)}{f_2(x)-f_2(y_1)}
\vert f_2(x)-f_2(y_1)\vert.
\]
It follows that $\vert\mu\vert=1$ and
\begin{gather*}
\vert(f_1+\mu f_2)(x)-(f_1+\mu f_2)(y_1)\vert=
\vert(f_1(x)-f_1(y_1))+\mu (f_2(x)-f_2(y_1))\vert=\\
\vert (f_1(x)-f_1(y_1))(1+\vert f_2(x)-f_2(y_1)\vert)\vert=
1+\vert f_2(x)-f_2(y_1)\vert >1,
\end{gather*}
which is untenable since the diameter of the range of
$f_1+\mu f_2$ is 1. This verifies our assertion.

\begin{step}\label{Sf}
Let $\{x_1,y_1\},\{x_2,y_2\}\in \X$. We have
$\{x_1,y_1\}\cap \{x_2,y_2\}=\emptyset$ if and only if
$G'(\{x_1,y_1\})\cap G'(\{x_2,y_2\})=\emptyset$.
\end{step}

The sufficiency follows from Step~\ref{Se}. As for the necessity,
observe that if $G'_{-1}$ denotes the function corresponding to
$\phi^{-1}$ like $G'$ corresponds to $\phi$ in
Step~\ref{Sc},
then we have $G'_{-1}={(G')}^{-1}$. Indeed, pick an $\{ x,y\}\in \X$.
Let $G'_{-1}(\{x,y\})=\{ a,b\}$ and $G'(\{ a,b\})=\{x',y'\}$. Applying
Step~\ref{Sd} for $\phi$ and then for $\phi^{-1}$, we find that for a
function $f\in C(X)$ with $\phi(f)\in \cS_s(\{ x',y'\})$ we have
$\phi^{-1}(\phi(f))=f\in \cS_s(\{ a,b\})$ and then
$\phi(f)\in \cS_s(\{x,y\})$. This results in $\{ x,y\}=\{ x',y'\}=
G'(G'_{-1}(\{x,y\}))$. The assertion is now obvious.

\begin{step}\label{Sg}
Let $x\in X$. There exists a unique element $g(x)\in X$ such that
$g(x)\in G'(\{x,y\})$ for every $x\neq y\in X$.
The function $g:X\to X$ is bijective and we have
$\{ g(x), g(y)\}=G'(\{ x,y\})$ $(\{ x,y\} \in \X)$.
\end{step}

Let $y_1,y_2\in X$ be such that $x\neq y_1,x\neq y_2$
and $y_1\neq y_2$. Let $g(x)$ denote the one and only element of the
set $G'(\{x,y_1\})\cap G'(\{x,y_2\})$ (see Step~\ref{Se}).
We are going to show that $g(x)\in G'(\{ x,y\})$ holds true for every
$x\neq y\in X$. If $X$ has only three elements, this is
obvious. Otherwise, pick
$x\neq y\in X$ for which $y\neq y_1$, $y\neq y_2$ and
suppose on the contrary that $g(x)\notin G'(\{ x,y\})$.
If $a_1,a_2\in X$ are so that
$G'(\{x,y_1\})=\{g(x),a_1\}$ and $G'(\{x,y_2\})=\{g(x),a_2\}$,
then taking the fact into account that by Step~\ref{Sf}
the sets $G'(\{x,y\})\cap G'(\{x,y_1\})$ and
$G'(\{x,y\})\cap G'(\{x,y_2\})$ are nonempty, we obtain that
$G'(\{x,y\})=\{a_1,a_2\}$. Similarly, we infer that
$G'(\{ y,y_1\})$ contains an element from $G'(\{ x,y\})=\{a_1,a_2\}$
as well as from $G'(\{ x,y_1\})=\{ g(x), a_1\}$. If
we had $g(x), a_2\in G'(\{y,y_1\})$, then we would obtain
$G'(\{x,y_2\})=\{ g(x), a_2\}=G'(\{y,y_1\})$ which would further imply
that
$\{x,y_2\}=\{y,y_1\}$, an apparent contradiction. Therefore, we have
$a_1\in G'(\{ y,y_1\})$ and a similar argument can be applied to
get $a_2\in G'(\{ y,y_2\})$. Hence, by Step~\ref{Sf} again, we
can choose a point $g(x)\neq b\in X$ so that
$G'(\{y,y_1\})=\{a_1,b\}$ and $G'(\{y,y_2\})=\{a_2,b\}$.
In the same fashion as above, one can check that
$G'(\{y_1,y_2\})=\{g(x),b\}$.
To sum up, we have the following relations
\begin{gather}
G'(\{ x,y_1\})=\{ g(x), a_1\}, \, G'(\{x,y_2\})=\{g(x),a_2\}, \,
G'(\{x,y\})=\{a_1,a_2\}, \notag \\
G'(\{y,y_1\})=\{a_1,b\}, \, G'(\{y,y_2\})=\{a_2,b\}, \,
G'(\{y_1,y_2\})=\{g(x),b\} \label{E:Sg1}.
\end{gather}
Assume for a moment that $X=\{x,y,y_1,y_2\}$. Let $f$ be the function
which takes
the following values: 0 at $y$, 1 at $y_1$, $e^{i\pi/3}$ at $y_2$ and
the geometric
center $(1+e^{i\pi/3})/3$ of the triangle $[0,1,e^{i\pi/3}]$ at $x$.
Let $\chi$ denote the characteristic function of the singleton $\{x\}$.
Since $\phi(\chi)$ is nonconstant, one of the values
$\phi(\chi)(a_1), \phi(\chi)(a_2), \phi(\chi)(g(x))$ differs from
$\phi(\chi)(b)$. Let it be $\phi(\chi)(g(x))$.
Since $f\in \cS(\{ y_1,y_2\})$ and $\diam(f(X))=1$, using
the last equation in the second line in \eqref{E:Sg1}, we have
$| \phi(f)(b)-\phi(f)(g(x))|=1$. Pick a nonzero $\mu\in
\mathbb{C}$ of modulus small enough to guarantee that $\mu +
(1+e^{i\pi/3})/3$ is still inside the triangle $[0,1, e^{i\pi/3}]$ and
for which
\[
\mu (\phi(\chi)(b)-\phi(\chi)(g(x)))=\lambda
(\phi(f)(b)-\phi(f)(g(x)))
\]
holds true with some positive scalar $\lambda$. Now, since
$f+\mu \chi\in \cS (\{ y_1,y_2\})$ and $\diam(f+\mu \chi)=1$,
by \eqref{E:Sg1} we have
\[
|\phi(f+\mu \chi)(b)-\phi(f+\mu \chi)(g(x))|=1.
\]
On the other hand, we can compute
\begin{multline*}
|\phi(f+\mu \chi)(b)-\phi(f+\mu \chi)(g(x))|=\\
|(\phi(f)(b)-\phi(g(x)))+\mu (\phi(\chi)(b)-\phi(
\chi)(g(x)))|=1+\lambda>1
\end{multline*}
which is a contradiction. If $X$ has at least five points, then by
Step~\ref{Se} we deduce that for
any point $z\in X, \, z\neq x,y,y_1,y_2$,
the set $G'(\{x,z\})$ has one common element with each of the sets
$G'(\{ x,y_1\})$, $G'(\{x,y_2\})$, $G'(\{x,y\})$. But having a look
at the first line in \eqref{E:Sg1}, one can see that this is a
contradiction again.
Thus, we have proved that $g(x)\in G'(\{x,y\})$ holds true for
every $x\neq y\in X$.

In view of Step~\ref{Se}, the uniqueness of $g(x)$ is obvious.

We next prove that $g$ is injective.
Suppose for a moment that $X$ has exactly three points $x,y,z$.
If, for example, $g(x)=g(y)=a$, then we have
\[
a\in G'(\{x,y\})\cap G'(\{x,z\})\cap G'(\{y,z\}).
\]
Now, using the facts that each of the sets
$G'(\{x,y\})$, $G'(\{x,z\})$, $G'(\{y,z\})$ has two different
elements, we easily arrive at a contradiction.
As for the case when $X$ has at least four points, we can argue as
follows. Let $x_1\neq x_2$ and assume that $g(x_1)=g(x_2)=b$. Choose
$\{y_1,y_2\}\in \X$ such that $\{x_1,x_2\} \cap \{y_1, y_2\}=\emptyset$.
We have $b\in G'(\{ x_1,y_1\})$ and $b\in G'(\{x_2,y_2\})$. But, on the
other hand, since the sets $\{ x_1,y_1\}$ and $\{ x_2,y_2\}$ are
disjoint, by Step~\ref{Sf} it follows that
$G'(\{ x_1,y_1\})\cap G'(\{x_2,y_2\})=\emptyset$. Since this is
a contradiction again, we have the injectivity of $g$.

We now show that $g$ is surjective. Let $x_0\in X$ be arbitrary
and let $y_0\in X,y_0\neq x_0$. By the surjectivity
of $G'$ there exists an $\{x,y\}\in \X$ for which
$\{g(x), g(y)\}=G'(\{x,y\})=\{x_0,y_0\}$. Therefore, $x_0$ is in the
range of $g$ verifying its surjectivity.

\begin{step}\label{Sh}
There exists a complex number $\tau$ of modulus 1 so that for every
$f\in \cT(x,y,u)$ we have $\phi(f)\in \cT(g(x),g(y),\tau u)$.
\end{step}

Let $\{x,y\}\in \X$ be arbitrary. By the surjectivity
of $\phi$ there exists an $f_0\in \cT(x,y,1)$ for which
$\phi(f_0)\in \cS_s(\{g(x),g(y)\})$.
Let
\[
\tau(x,y)=\phi(f_0)(g(x))-\phi(f_0)(g(y)).
\]
Then we have
\[
H(x,y,1)\subseteq T(\phi(f_0))=\{(g(x),g(y),\tau(x,y)),
(g(y),g(x),-\tau(x,y))\}.
\]
By the definition of $H$ and its nonemptyness (Step~\ref{Sa}) we obtain
\[
H(x,y,1)=\{(g(x),g(y),\tau(x,y)),
(g(y),g(x),-\tau(x,y))\}.
\]
It is now easy to see that the implication
\begin{equation}\label{E:Sh1}
f\in \cT(x,y,u) \Longrightarrow
\phi(f)\in \cT(g(x),g(y),\tau(x,y) u)
\end{equation}
holds true for every $u\in \mathbb{C}$.
It remains to verify that $\tau$ does not depend on its
variables $x,y$.
To this end first observe that, by \eqref{E:Sh1}, $\tau$ is symmetric,
i.e. $\tau(x,y)=\tau(y,x)$ $(\{ x,y\}\in \X)$.
Next, let $x,y_1,y_2\in X$ be pairwise different points.
Pick functions $f_1\in \cT(x,y_1,-1)$ and $f_2\in \cT(x,y_2,-1)$
with disjoint supports and with ranges in [0,1].
Define
\[
f=f_1+e^{i\pi/3} f_2\in \cT(x,y_1,-1)
\cap \cT(x,y_2,-e^{i\pi/3})\cap \cT(y_1,y_2, 1-e^{i\pi/3}).
\]
Then we have
\begin{gather*}
\phi(f)\in \cT(g(x),g(y_1),-\tau(x,y_1))\cap \\
\cT(g(x),g(y_2),-e^{i\pi/3} \tau(x,y_2))\cap
\cT(g(y_1), g(y_2), \tau(y_1,y_2)(1-e^{i\pi/3})).
\end{gather*}
Hence, for the complex numbers
$\phi(f)(g(x))$, $\phi(f)(g(y_1))$ and $\phi(f)(g(y_2))$ we have
\begin{gather*}
|\phi(f)(g(x))-\phi(f)(g(y_1))|=|\phi(f)(g(x))-\phi(f)(g(y_2))|=\\
|\phi(f)(g(y_1))-\phi(f)(g(y_2))|=1.
\end{gather*}
It is easy to see that for any complex numbers
$a,b$ of modulus 1, if $|a-b|=1$, then we have $a=e^{\pm i\pi/3}b$.
Therefore, we obtain that
\begin{gather*}
-\tau(x,y_1)=
\phi(f)(g(x))-\phi(f)(g(y_1))=\\
e^{\pm i\pi/3}(\phi(f)(g(x))-\phi(f)(g(y_2)))=
e^{\pm i\pi/3}(-e^{i\pi/3}\tau(x,y_2)),
\end{gather*}
i.e. $\tau(x,y_1)=e^{\pm i\pi/3}e^{i\pi/3}\tau(x,y_2)$.
Applying the same argument to the function $f_1+e^{-i\pi/3}f_2$, we
obtain $\tau(x,y_1)=e^{\pm i\pi/3}e^{-i\pi/3}\tau(x,y_2)$.
Comparing these two equalities, we conclude
$\tau(x,y_1)=\tau(x,y_2)$. Taking the symmetry of $\tau$ into account,
it is now apparent that $\tau$ is a constant function. Let
this constant be denoted by the same symbol $\tau$.
Our assertion follows from \eqref{E:Sh1}.

\begin{step}\label{Si}
For every $f\in C(X)$, the function
$\phi(f\circ g)-\tau \cdot f$ is constant.
\end{step}

Let $f\in \cT(x,y,1)$. Then by Step~\ref{Sh} it follows that
$\phi(f)\in \cT(g(x),g(y),\tau)$. Thus we have
\[
\phi(f)(g(x))-\phi(f)(g(y))=\tau=
\tau (f(x)-f(y)),
\]
which implies
\begin{equation}\label{E:Si1}
\phi(f)(g(y))-\tau f(y)= \phi(f)(g(x))-\tau f(x).
\end{equation}
Let $z\in X$ be so that $z\neq x,y$.
Set $u=f(x)-f(z)$. Clearly, $|u|\leq 1$.
If $u=0$, then $f\in \cT (z,y,1)$ which further implies
$\phi(f)\in \cT(g(z),g(y),\tau)$. Just as on the way to \eqref{E:Si1},
we obtain that
\[
\phi(f)(g(z))-\tau f(z)=
\phi(f)(g(y))-\tau f(y)=
\phi(f)(g(x))-\tau f(x).
\]
Suppose now that $u\neq 0$. Define
\[
U=\left\{ p\in X\setminus\{y\}\, :\, f(p)\neq f(x)
\text{ and }
\left\vert \frac{f(x)-f(p)}{\vert f(x)-f(p)\vert}-
\frac{u}{\vert u \vert}\right\vert < \frac{1}{2}\right\}.
\]
Since $f$ is continuous, $U$ is an open neighbourhood of the point $z$.
By Uryson's lemma, there exists a function
$f_0\in C(X)$ with range in $[0,1]$ and support in $U$ for which
$f_0(z)=1$ and $f_0(y)=0$. Let
\[
f_1=\frac{\vert u\vert f_0}
{\max\{ \vert f(x)-f\vert,\vert u\vert\} }
\, \text{ and }\,
f_2=f_1(f(x)-f).
\]
Clearly, the support of $f_2$ is
included in that of $f_0$ which is a subset of $U$.
By the definition of $U$ and $f_2$ it is now quite easy to
check that $\diam (f_2(X))\leq |u|$. On the other hand,
$f_2(z)-f_2(x)=u-0$ and $f_2(z)-f_2(y)=u-0$. Consequently, we have
\begin{equation}\label{E:Si2}
f_2\in \cT(z,x,u)\cap \cT(z,y,u).
\end{equation}
Since the range of $f_1$ is a subset of $[0,1]$, for an
arbitrary $w\in X$ we have
\[
(f+f_2)(w)=f(w)+f_1(w)(f(x)-f(w))=
(1-f_1(w))f(w)+f_1(w)f(x).
\]
This shows that $(f+f_2)(w)$ belongs to the convex hull of $f(X)$.
The diameter of this convex hull is equal to $\diam(f(X))$. Therefore,
we have
\[
\diam (f+f_2)(X))\leq \diam (f(X))=1.
\]
On the other hand, we infer
\begin{gather*}
(f+f_2)(z)-(f+f_2)(y)=f(z)-f(y)+f_2(z)-f_2(y)=\\
f(z)-f(y)+u= f(z)-f(y)+f(x)-f(z)=f(x)-f(y)=1.
\end{gather*}
Consequently,
\begin{equation}\label{E:Si3}
f+f_2\in \cT(z,y,1).
\end{equation}
By \eqref{E:Si2} and \eqref{E:Si3} we obtain that
\[
\phi(f_2)\in \cT(g(z),g(x),\tau u)\cap
\cT(g(z),g(y),\tau u),\, \phi(f+f_2)\in \cT(g(z),g(y),\tau).
\]
Hence, we can compute
\begin{gather*}
\phi(f)(g(x))-\phi(f)(g(z))=\\
(\phi(f)(g(x))-\phi(f)(g(y)))- (\phi(f)(g(z))-\phi(f)(g(y)))=\\
\tau-(\phi (f+f_2)(g(z))-\phi (f+f_2)(g(y)))+
(\phi(f_2)(g(z))-\phi(f_2)(g(y)))=\\
\tau-\tau +\tau u=\tau (f(x)-f(z)).
\end{gather*}
Thus, we have proved that
\[
\phi(f)(g(z))-\tau f(z)= \phi(f)(g(x))-\tau f(x)
\]
holds true for every $z\in X$.

Now, the proof of the theorem can be completed as follows.
By the linearity of $\phi$, there is a linear functional $t':C(X)\to \C$
such that
\[
\phi(f\circ g)-\tau \cdot f=t'(f)1 \qquad (f\in C(X)).
\]
Since $g$ is a bijection, with the notation $\varphi=g^{-1}$ and
$t(f)=t'(f\circ \varphi)$ we have
\begin{equation}\label{E:Si4}
\phi(f)-\tau \cdot f\circ \varphi=t(f)1 \qquad (f\in C(X)).
\end{equation}
It follows from \eqref{E:Si4}
that $f\circ \varphi$ is continuous for every $f\in C(X)$. Using
Uryson's lemma, we deduce that
$\varphi$ is continuous and hence, being a bijection between compact
Hausdorff spaces, we obtain that it is a homeomorphism.
Since the relation $t(1)\neq -\tau$ is obvious, the proof of our theorem
is complete.
\end{proof}

\begin{remark}
We believe that it would be interesting to study our problem for the
operator algebra $B(H)$ of all bounded linear operators acting on a
complex
Hilbert space $H$. We mean to characterize those linear bijections of
$B(H)$ which preserve the diameter of the spectrum.
The spectrum preserving, respectively spectral radius preserving
linear bijections of $B(H)$ were characterized by Jafarian and Sourour
\cite {JaSo}, respectively by Bre\v sar and \v Semrl \cite{BrSe1}.
\end{remark}


\end{document}